\documentclass[twoside,leqno,10pt]{amsart}
\usepackage{amsfonts}
\usepackage{amsmath}
\usepackage{amscd}
\usepackage{amssymb}
\usepackage{amsthm}
\usepackage{amsrefs}
\usepackage{latexsym}
\usepackage{mathrsfs}
\usepackage{bbm}
\usepackage{enumerate}
\usepackage{color}
\begin{document}

\newtheorem{theorem}[subsection]{Theorem}
\newtheorem{proposition}[subsection]{Proposition}
\newtheorem{lemma}[subsection]{Lemma}
\newtheorem{corollary}[subsection]{Corollary}
\newtheorem{conjecture}[subsection]{Conjecture}
\newtheorem{prop}[subsection]{Proposition}
\numberwithin{equation}{section}
\newcommand{\mr}{\ensuremath{\mathbb R}}
\newcommand{\mc}{\ensuremath{\mathbb C}}
\newcommand{\dif}{\mathrm{d}}
\newcommand{\intz}{\mathbb{Z}}
\newcommand{\ratq}{\mathbb{Q}}
\newcommand{\natn}{\mathbb{N}}
\newcommand{\comc}{\mathbb{C}}
\newcommand{\rear}{\mathbb{R}}
\newcommand{\prip}{\mathbb{P}}
\newcommand{\uph}{\mathbb{H}}
\newcommand{\fief}{\mathbb{F}}
\newcommand{\majorarc}{\mathfrak{M}}
\newcommand{\minorarc}{\mathfrak{m}}
\newcommand{\sings}{\mathfrak{S}}
\newcommand{\fA}{\ensuremath{\mathfrak A}}
\newcommand{\mn}{\ensuremath{\mathbb N}}
\newcommand{\mq}{\ensuremath{\mathbb Q}}
\newcommand{\half}{\tfrac{1}{2}}
\newcommand{\f}{f\times \chi}
\newcommand{\summ}{\mathop{{\sum}^{\star}}}
\newcommand{\chiq}{\chi \bmod q}
\newcommand{\chidb}{\chi \bmod db}
\newcommand{\chid}{\chi \bmod d}
\newcommand{\sym}{\text{sym}^2}
\newcommand{\hhalf}{\tfrac{1}{2}}
\newcommand{\sumstar}{\sideset{}{^*}\sum}
\newcommand{\sumprime}{\sideset{}{'}\sum}
\newcommand{\sumprimeprime}{\sideset{}{''}\sum}
\newcommand{\shortmod}{\ensuremath{\negthickspace \negthickspace \negthickspace \pmod}}
\newcommand{\V}{V\left(\frac{nm}{q^2}\right)}
\newcommand{\sumi}{\mathop{{\sum}^{\dagger}}}
\newcommand{\mz}{\ensuremath{\mathbb Z}}
\newcommand{\leg}[2]{\left(\frac{#1}{#2}\right)}
\newcommand{\muK}{\mu_{\omega}}

\title[Mean values of cubic and quartic {D}irichlet characters]{Mean values of cubic and quartic {D}irichlet characters}

\date{\today}
\author{Peng Gao and Liangyi Zhao}

\begin{abstract}
We evaluate the average of cubic and quartic Dirichlet character sums with the modulus going up to a size comparable to the length of the individual sums. This generalizes a result of Conrey, Farmer and Soundararajan \cite{CFS} on quadratic Dirichlet character sums.
\end{abstract}

\maketitle


\section{Introduction}

 Estimations for character sums have wide applications in number theory. A fundamental result is the well-known P\'olya-Vinogradov inequality (see for example \cite[Chap. 23]{Da}), which asserts that for any non-principal Dirichlet character $\chi$ modulo $q$, $M \in \mz$ and $N \in \mn$,
\begin{align} \label{pv}
   \sum_{M < n \leq M+N} \chi(n) \ll q^{1/2}\log q.
\end{align}

  One may regard the P\'olya-Vinogradov inequality as a mean value estimate for characters and one expects to obtain asymptotic formulas if an extra average over the characters is introduced. We are thus led to the investigation on the following expression
\begin{align*}
   \sum_{\chi \in S}\sum_{n} \chi(n)
\end{align*}
  where $S$ is a certain set of characters. For example, when $S$ is the set of all Dirichlet characters modulo $q$, then it follows from the orthogonality relation of the characters that the above sum really amounts to a counting on the number of integers which are congruent to $1$ modulo $q$. \newline

  Another interesting choice for $S$ is to take $S$ to contain characters of a fixed order. A basic and important case is the set of quadratic characters.
In \cite{CFS}, J. B. Conrey, D. W. Farmer and K. Soundararajan studied the following sum:
\begin{align*}
   S_2(X,Y)= \sum_{\substack {m \leq X \\ (m, 2)=1}}\sum_{\substack {n \leq Y \\ (n, 2)=1}} \leg {m}{n},
\end{align*}
  where $\leg {m}{n}$ is the Jacobi symbol. \newline

   While it is relatively easy to obtain an asymptotic formula of $S_2(X, Y)$ when $Y=o(X/\log X)$ or $X=o(Y/\log Y)$ using \eqref{pv}, it is more subtle to treat the remaining $XY$-ranges. Using a Poisson summation formula developed in \cite{sound1}, a valid asymptotic formula of $S_2(X, Y)$ for all $X,Y$ is obtained in \cite{CFS}.  Most interestingly, the formula exhibits a transition in the behavior of $S_2(X,Y)$ when $X$ and $Y$ are of comparable size. \newline

  Recently, the authors \cite{G&Zhao2019} studied the mean values of some quadratic, cubic and quartic {H}ecke characters. These are analogues to $S_2(X,Y)$  in number fields. Again, the most interesting case is when $X, Y$ are of comparable size. Another similar
behavior is shown to exist by I. Petrow in an earlier study \cite{Petrow} on the mean values of shifted convolution sums of Fourier coefficients of Hecke eigenforms. \newline

   In this paper, we return to the classical setting by studying the mean values of cubic and quartic Dirichlet characters.
   To form such sums, we note that primitive cubic Dirichlet characters exist modulo any rational prime $p$ if and only if $p \equiv 1 \pmod 3$ and
    there are precisely two such characters which are complex conjugate of each other when they exist.  The same conclusions apply to primitive quartic Dirichlet characters as well except in we need to have $p \equiv 1 \pmod 4$.   Based on these observations, we introduce the following sets of Dirichlet characters that we aim to study. We define the set $S_{3,1}$ to be the set that contains the principal Dirichlet character modulo $1$. For any integer $n > 1$, we define the set $S_{3,n}$ to be non-empty if and only if $n$ is a product of powers of primes which are congruent to $1$ modulo $3$, in which case by writing $n=\prod^{k}_{i=1}p^{\alpha_i}_i$ with $p_i \equiv 1 \pmod 3$ and $\alpha_i \geq 1$, we define
\begin{align} \label{S3}
   S_{3,n}=\left\{ \prod^k_{i=1}\chi^{\alpha_i}_i: \chi_i \hspace{0.05in} \text{primitive cubic Dirichlet character modulo $p_i$} \right\}.
\end{align}

   Similarly, let $S_{4,1}$ be the set that contains the principal Dirichlet character modulo $1$. For any integer $n > 1$, we define the set $S_{4,n}$ to be non-empty if and only if $n$ is a product of powers of primes which are congruent to $1$ modulo $4$, in which case by writing $n=\prod^{k}_{i=1}p^{\alpha_i}_i$ with $p_i \equiv 1 \pmod 4$ and $\alpha_i \geq 1$, we define
\begin{align} \label{S4}
   S_{4,n}=\left\{ \prod^k_{i=1}\chi^{\alpha_i}_i: \chi_i \hspace{0.05in} \text{primitive quartic Dirichlet character modulo $p_i$} \right\}.
\end{align}

   We are now ready to define the following character sums of our interest. For $i=3$, $4$, $X$, $Y>1$, let
\begin{align*}
  S_i(X,Y) =\sum_{n \leq Y}\sum_{\chi \in S_{i,n}} \sum_{m \leq X} \chi(m).
\end{align*}

   Our goal in this paper is to evaluate $S_i(X,Y)$ asymptotically for $i=3,4$. Our method here only allows us to treat the situation in which $Y \leq X$, a condition we shall assume henceforth.  In this case, one expects that the main contribution for $S_i(X,Y)$ comes from the terms when $n$ is a cube if $i=3$ or a fourth-power if $i=4$. Treating the remaining terms using the P\'olya-Vinogradov inequality \eqref{pv}, we deduce that
\begin{align}
\label{Sieleasmp}
  S_i(X,Y) = C_i\frac {XY^{1/i}}{\sqrt{\log Y}}+O \left( Y^{3/2+\epsilon} \right) , \quad i=3,4,
\end{align}
   for some constants $C_i$. The above allows us to obtain asymptotic formulas for $S_i(X, Y), i=3,4$ when $Y$ is small compared to $X$.
   Thus, it is again a subtlety to obtain asymptotic formulas for $S_i(X, Y), i=3,4$ when $Y$ is close to $X$.
   This is precisely what we want to study in this paper. In view of \eqref{Sieleasmp},
   we may assume that $Y \geq X^{6/7}$ when studying $S_3(X,Y)$ and that $Y \geq X^{4/5}$ when studying $S_4(X,Y)$.   Our main result is
\begin{theorem}
\label{cubicquarticmean}
   For large $X \geq Y$ and any $\varepsilon>0$, we have for $Y \geq X^{6/7}$,
\begin{align}
\label{S3asymp}
   S_3(X,Y)= C_1\frac {XY^{1/3}}{\sqrt{\log Y}}+O \left (\frac {XY^{1/3}}{(\log Y)^{3/2}}+Y^{4/3}+
   Y^{24/17+\epsilon}\left ( \frac {Y}{X} \right )^{5/17} \right ).
\end{align}
   We also have for $Y \geq X^{4/5}$,
\begin{align}
\label{S4asymp}
   S_4(X,Y)= C_2\frac {XY^{1/4}}{\sqrt{\log Y}}+O \left (\frac {XY^{1/4}}{(\log Y)^{3/2}}+Y^{91/62+\epsilon}\left(\frac {Y}{X} \right)^{9/31}+
   X Y^{13/31+\epsilon} \left ( \frac {Y}{X} \right )^{38/31} \right ) .
\end{align}
Here $C_1, C_2$ are constants given in \eqref{C1} and \eqref{C2}, respectively.
\end{theorem}

A simple computation will show that \eqref{S3asymp} is a valid asymptotic formula (the main term dominates all the $O$-terms) if $Y = O(X^{33/35-\delta})$ and \eqref{S4asymp} if $Y =O( X^{160/187-\delta})$ for any $\delta >0$. \newline

We point out here that $C_1$ and $C_2$ are positive (see the discussions below \eqref{C1} and \eqref{C2}). Our proof of Theorem~\ref{cubicquarticmean} follows the line of treatment in \cite{CFS} by first applying the Poisson summation formula to covert the sum over $m$ in $S_i(X,Y)$ to its dual sum. The zero frequency gives us the main term as usual. To treat the contribution of remaining terms, we transform the sum over $n$ in $S_i(X,Y)$, using Lemma \ref{lemma:cubicquarticclass}, into another sum over algebraic integers in an imaginary quadratic number field. The resulting sums involve with certain Gauss sums for which we shall use a result of S. J. Patterson \cite{P} (see Lemma~\ref{lem1}). This treatment is also inspired by the method used in \cite{B&Y}.

\subsection{Notations} The following notations and conventions are used throughout the paper.\\
\noindent $e(z) = \exp (2 \pi i z) = e^{2 \pi i z}$, $\omega=e(1/3)$. \newline
$f =O(g)$ or $f \ll g$ means $|f| \leq cg$ for some unspecified
positive constant $c$. \newline
$\varepsilon$ denotes an arbitrary small positive number, which may be different
from line to line. \newline
$N(n)$ denotes the norm of $n \in \mz[\omega]$ or $n \in \mz[i]$. \newline
$\varpi$ denotes a prime element in $\mz[\omega]$ or $\mz[i]$ and $p$ denotes a prime in $\mz$. \newline

\section{Preliminaries}
\label{sec 2}

\subsection{Cubic and quartic symbols}
\label{sec2.4}
    For any number field $K$, let $\mathcal{O}_K$ denote the ring of integers in $K$ and $U_K$ the group of units in $\mathcal{O}_K$. Throughout this paper, set $K_{\omega}=\mq(\omega), K_i=\mq(i)$, where $\omega=\exp(2\pi i/3)$. It is well-known that both $\mq(i)$ and $\mq(\omega)$ have class number one and that $\mathcal{O}_{K_{\omega}}=\mz[\omega], \mathcal{O}_{K_i}=\mz[i]$. Recall that every ideal in $\intz[\omega]$ co-prime to $3$ has a unique generator congruent to $1$ modulo $3$ (see \cite[Proposition 8.1.4]{BEW}) and every ideal in $\intz[i]$ coprime to $2$ has a unique generator congruent to $1$ modulo $(1+i)^3$ (see the paragraph above Lemma 8.2.1 in \cite{BEW})). These generators are called primary. \newline

Let $\leg{\cdot}{n}_3$ be the cubic residue symbol in $\mathcal{O}_{K_{\omega}}$.  For a prime $\varpi \in \mathcal{O}_{K_{\omega}}$
with $N(\varpi) \neq 3$, the cubic symbol is defined for $a \in \mathcal{O}_{K_{\omega}}$ , $(a, \varpi)=1$ by $\leg{a}{\varpi}_3 \equiv
a^{(N(\varpi)-1)/3} \pmod{\varpi}$, with $\leg{a}{\varpi}_3 \in \{
1, \omega, \omega^2 \}$. When $\varpi | a$, we define
$\leg{a}{\varpi}_3 =0$.  Then the cubic symbol can be extended
to any composite $n$ with $(N(n), 3)=1$ multiplicatively. We extend the definition of $\leg{\cdot }{n}_3$ to $n \in U_{K_{\omega}}$ by setting $\leg{\cdot}{n}_3=1$. Recall that \cite[Theorem 7.8]{Lemmermeyer} the cubic reciprocity law states
that for two primary $m, n \in \mathcal{O}_{K_{\omega}}$,
\begin{equation*}
 \leg{m}{n}_3 = \leg{n}{m}_3.
\end{equation*}

The quartic case is similar.   Let $\leg {\cdot}{n}_4$ be the quartic residue symbol in $\mathcal{O}_{K_{i}}$.  Suppose that $\varpi \in \mathcal{O}_{K_{i}}$ is a prime with $N(\varpi) \neq 2$.   If $a \in
\mathcal{O}_{K_{i}}$, $(a, \varpi)=1$, then we define $\leg{a}{\varpi}_4$ by $ \leg{a}{\varpi}_4 \equiv
a^{(N(\varpi)-1)/4} \pmod{\varpi}$, with $\leg{a}{\varpi}_4 \in \{ \pm 1, \pm i \}$. Set $\leg{a}{\varpi}_4 =0$ if $\varpi | a$.  Now the quartic character can be extended to any composite $n$ with $(N(n), 2)=1$ multiplicatively.  As before, we extend the definition of $\leg{\cdot }{n}_4$ to $n \in U_{K_{i}}$ by setting $\leg{\cdot}{n}_4=1$. The quartic reciprocity law, \cite[Theorem 6.9]{Lemmermeyer}, states that for two primary $m, n \in \mathcal{O}_{K_{i}}$,
\begin{equation} \label{quartrec}
 \leg{m}{n}_4 = \leg{n}{m}_4(-1)^{((N(n)-1)/4)((N(m)-1)/4)}.
\end{equation}

   Similar to \cite[Lemma 2.1]{B&Y} and \cite[Section 2.2]{G&Zhao}, we have the following description of $S_{i,n}$, $i=3$, $4$ defined in \eqref{S3} and \eqref{S4} using the cubic and quartic residue symbols.
\begin{lemma}
\label{lemma:cubicquarticclass}
  For any non-empty $S_{3,n}$, there is a bijection between $S_{3,n}$ and the set of cubic residue symbols of the form $\chi_{3,q}:m \rightarrow (\frac{m}{q})_3$ for some $q \in \mathcal{O}_{K_{\omega}}$ with $q$ primary, not divisible by any $\mq$-rational primes and $N(q) = n $.
  Similarly, a bijection exists between any non-empty $S_{4,n}$ and the set of quartic residue symbols of the form $\chi_{4,q}:m \rightarrow (\frac{m}{q})_4$ for some $q \in \mathcal{O}_{K_{i}}$ with $q$ primary, not divisible by any $\mq$-rational primes and $N(q) = n $.
\end{lemma}

\subsection{Gauss sums}
\label{section:Gauss}

  Let $\chi$ be a Dirichlet character of modulus $n$.  For any $r \in \mz$, we define the Gauss sum $\tau(r, \chi)$ as follows:
\begin{align}
\label{taur}
  \tau(r, \chi)=\sum_{x \bmod {n}}\chi(x)e\left( \frac{rx}{n} \right).
\end{align}

  Similarly, for any $n, r \in \mathcal{O}_{K_{\omega}}$, we define
\begin{align*}
 g_3(r,n) = \sum_{x \bmod{n}} \leg{x}{n}_3 \widetilde{e}_{\omega}\leg{rx}{n}, \quad \mbox{where} \quad \widetilde{e}_{\omega}(z) =\exp \left( 2\pi i  \left( \frac {z}{\sqrt{-3}} - \frac {\bar{z}}{\sqrt{-3}} \right) \right).
\end{align*}
Furthermore, for any $n, r \in \mathcal{O}_{K_{i}}$, set
\begin{align*}
  g_4(r,n) = \sum_{x \bmod{n}} \leg{x}{n}_4 \widetilde{e}_i\leg{rx}{n}, \quad \mbox{with} \quad \widetilde{e}_i(z) =\exp \left( 2\pi i  \left( \frac {z}{2i} - \frac {\bar{z}}{2i} \right) \right) .
\end{align*}

The following (or similar) properties of $g_i(r,n)$ for $i=3,4$ that we enumerate below can be found in \cite[Section 3]{H&P} for $i=3$ and \cite[Section 2]{Diac} for $i=4$.  Our $g_3(r,n)$ differs slightly from the cubic Gauss sums defined in \cite{H&P}.  But the properties of $g_3(r,n)$ can be proved in the same way as those of its analogue in \cite{H&P}.  Indeed, both \eqref{eq:gmult} for $i=3$ and \eqref{gprod} can be proved using the cubic reciprocity law and the Chinese remainder theorem.

\begin{align}
\label{eq:gmult}
 g_i(rs,n) & = \overline{\leg{s}{n}}_i g_i(r,n), \quad (s,n)=1, \qquad \mbox{$n$ primary}.
\end{align}

Note that we have \cite[(13)]{B&Y} for $(n_1, n_2)=1$ and $n_1, n_2$ primary,
\begin{align} \label{gprod}
   g_3(k, n_1n_2)=\overline{\leg {n_1}{n_2}}_3 g_3(k, n_1)g_3(k,n_2).
\end{align}
Similarly, \cite[Lemma 2.3]{G&Zhao4} gives that for $(n_1, n_2)=1$ and $n_1, n_2$ primary,
\begin{align}
\label{g4prod}
    g_4(r,n_1 n_2) =\leg{n_2}{n_1}_4\leg{n_1}{n_2}_4g_4(r, n_1) g_4(r, n_2)=(-1)^{((N(n_1)-1)/4)((N(n_2)-1)/4)}\leg{n^2_1}{n_2}_4g_4(r, n_1) g_4(r, n_2),
\end{align}
  where the last equality above follows from the quartic reciprocity law, \eqref{quartrec}. \newline

It is well-known (see \cite[(11)]{B&Y} and \cite[Prop. 6.5]{Lemmermeyer}) that for a primary prime $\varpi$ (belonging to the corresponding ring of integers),
\begin{align*}
    |g_i(1,\varpi)|=N(\varpi)^{1/2}, \quad i=3,4.
\end{align*}
    Using \eqref{gprod} or \eqref{g4prod} to reduce the general case to the prime case and applying the above bound, we deduce that
    when $r, n$ are in the corresponding ring of integers with $n$ being square-free,  $g_i(r,n) \neq 0$, $i=3$, $4$ only when $(r,n)=1$, in which case we get
\begin{align}
\label{grnbound}
    |g_i(r,n)| \leq N(n)^{1/2}, \quad i=3,4.
\end{align}

  In the proof of Theorem \ref{cubicquarticmean}, we need a relation between $\tau$ and $g_i$, $i=3$, $4$. Analogous to the discussions in \cite[Section 2.2]{B&Y} (but be aware that the notion of $g_3(r,n)$ there is slightly different from ours), one shows that for $r \in \mz$ and any $\chi \in S_{3,n}$ which corresponds to $\chi_{3,q}$,
\[  \tau(r, \chi)=\sum_{x \bmod{N(q)}} \leg{x}{q}_3 e \left( \frac{rx}{N(q)} \right) . \]
Now we can write $x=y \bar{q} + \bar{y} q \pmod{q \bar{q}}$ with $y$ running through a set of representative $\pmod{q}$ in $\intz[\omega]$.  Here $\bar{q}$ is the complex conjugate of $q$.  Every $x \pmod{N(q)}$ is uniquely representable this way.  So we easily get
\[ \tau(r, \chi) = \sum_{y \bmod{q}} \leg{y\bar{q}}{q}_3 e \left( r \left( \frac{y}{q} + \frac{\bar{y}}{q} \right) \right) = \sum_{y \bmod{q}} \leg{y}{q}_3 e \left( r \left( \frac{y}{q} + \frac{\bar{y}}{q} \right) \right)  , \]
upon noting that $\leg{\bar{q}}{q}_3 =1$ by cubic reciprocity.  Now the change of variables $y \to y \sqrt{-3}$ gives
\begin{align}
\label{taug3}
  \tau(r, \chi) = \overline{\leg {\sqrt{-3}}{q}}_3 g_3(r, q).
\end{align}

A similar relation exists for quartic Gauss sums.  For $r \in \mz$ and any $\chi \in S_{4,n}$ which corresponds to $\chi_{4,q}$ by Lemma \ref{lemma:cubicquarticclass}, we have (\cite[p. 894]{G&Zhao})
\begin{equation}
\label{tau}
\tau(r, \chi)=\overline{\leg {2i}{q}}_4 \leg {\overline{q}}{q}_4 g_4(r, q).
\end{equation}

  In order to apply \eqref{tau}, we need to express $\leg {\overline{q}}{q}_4$ in terms of ray class characters $\pmod {16}$ in $\mq(i)$.
  To that end, note that we have $q \equiv 1 \pmod {(1+i)^3}$ with $q$ having no rational prime divisors. If we write
  $q=a+bi$ with $a,b \in \mz$, then we deduce that $(a, b)=1$ so that
\begin{align*}
  \leg {\overline{q}}{q}_4= \leg {a-bi}{a+bi}_4=\leg {2a}{a+bi}_4=\leg {2(-1)^{(N(n)-1)/4}}{a+bi}_4\leg {(-1)^{(N(n)-1)/4}a}{a+bi}_4.
\end{align*}

   We further observe that
$q=a+bi$, $a$, $b \in \mz $ in $\mz[i]$ is congruent to $1 \bmod{(1+i)^3}$ if and only if $a \equiv 1 \pmod{4}$, $b \equiv 0 \pmod{4}$ or $a \equiv 3 \pmod{4}$, $b \equiv 2 \pmod{4}$ by \cite[Lemma 6, p. 121]{I&R}. It follows from this that we have
\begin{align} \label{a&b}
  a \equiv (-1)^{(N(q)-1)/4} \pmod 4,  \quad b \equiv 1-(-1)^{(N(q)-1)/4} \pmod 4.
\end{align}

   As $(-1)^{(N(q)-1)/4}a$ is primary according to \eqref{a&b}, we have by the quartic reciprocity law,
\begin{align*}
   \leg {(-1)^{(N(q)-1)/4}a}{a+bi}_4=(-1)^{((N(a)-1)/4)((N(q)-1)/4)} \leg {a+bi}{a}_4= \leg {a+bi}{a}_4= \leg {bi}{a}_4= \leg {b}{a}_4 \leg {i}{a}_4= \leg {i}{a}_4,
\end{align*}
   where the last equality follows from \cite[Proposition 9.8.5]{I&R}, which states that for $a, b \in \mz, (a, 2b)=1$ ,
\begin{align*}
   \leg {b}{a}_4=1.
\end{align*}

One of the supplement laws to the quartic reciprocity law states that if $n=a+bi$ being primary, then
\begin{align}
\label{2.05}
  \leg {i}{n}_4=i^{(1-a)/2}.
\end{align}

   It follows from \eqref{2.05} that
\begin{align*}
   \leg {i}{a}_4=i^{(1-(-1)^{\frac {N(q)-1}{4}}a)/2}=(-1)^{(a^2-1)/8} =: \lambda_0(q).
\end{align*}
It is easy to check that $\lambda_0$ is a ray class character $\pmod {16}$ in $\mq(i)$. \newline

   On the other hand, note that by the definition that
\begin{align*}
   \leg {(-1)^{(N(q)-1)/4}}{a+bi}_4=(-1)^{\frac {N(q)-1}{4}\cdot \frac {N(q)-1}{4}}=(-1)^{(N(q)-1)/4}=\leg {-1}{q}_4.
\end{align*}

   We then deduce that
\begin{align*}
  \leg {\overline{q}}{q}_4= \leg {-2}{q}_4\lambda_0(q).
\end{align*}

   We conclude from this and \eqref{tau} that for $r \in \mz$ and any $\chi \in S_{4,n}$ which corresponds to $\chi_{4,q}$ by
Lemma \ref{lemma:cubicquarticclass}, we have
\begin{align*}
  \tau(r, \chi)=\overline{\leg {-i}{q}}_4 \lambda_0(q)g_4(r, q).
\end{align*}

\subsection{Analytic behavior of Dirichlet series associated with Gauss sums}
\label{section: smooth Gauss}

    In the proof of Theorem \ref{cubicquarticmean}, we need to know the analytic behavior of certain Dirichlet series associated with cubic or quartic Gauss sums.
    For any  ray class character $\psi_3$
   $\pmod {9}$ in $\mq(\omega)$ and any ray class character $\psi_4$ $\pmod {16}$ in $\mq(i)$, we define
\begin{align*}
G_3(s,dk; \psi_3)=\sum_{\substack{n \equiv 1 \bmod {3}\\ (n,d)=1}} \frac { \psi_3(n) g_3(dk,n)}{N(n)^{s}}, \quad
G_4(s,dk;\psi_4) =\sum_{\substack{n \equiv 1 \bmod {(1+i)^3} \\ (n,d)=1}}  \frac {\psi_4(n) g_4(dk,n)}{N(n)^{s}}.
\end{align*}

  We deduce from a general result of S. J. Patterson \cite[Lemma, p. 200]{P} the following analytic behavior of $G_i$, with $i=3,4$ (see also \cite[Lemma 3.5]{B&Y}).
\begin{lemma}
\label{lem1} The functions $G_i(s,dk; \psi_i), i=3,4$ has meromorphic continuation to the half plane with $\Re (s) > 1$.  It is holomorphic in the
region $\sigma=\Re(s) > 1$ except possibly for a pole at $s = 1+1/i$. For any $\varepsilon>0$, letting $\sigma_1 = 3/2+\varepsilon$,
then for $\sigma_1 \geq \sigma \geq \sigma_1-1/2$, $|s-(1+1/i)|>1/(2i)$, we have
\[ G_i(s,dk;\psi_i) \ll N(dk)^{\frac 12(\sigma_1-\sigma+\varepsilon)}(1+t^2)^{\frac {i-1}2(\sigma_1-\sigma+\varepsilon)}, \]
  where $t=\Im(s)$ and the norm is taken in the corresponding number field. Moreover, the residue satisfies
\[ \mathrm{Res}_{s=1+1/3}G_3(s,dk;\psi_3) \ll N((dk)_1)^{-1/6+\varepsilon}, \quad \mathrm{Res}_{s=1+1/4}G_4(s,dk;\psi_4) \ll N(dk)^{1/8+\varepsilon}, \]
 where we write $dk=(dk)_1(dk)^2_2(dk)^3_3$ with $(dk)_1(dk)^2_2$ cubic-free in $\mz[\omega]$.
\end{lemma}

\section{Proof of Theorem \ref{cubicquarticmean}}
\label{sec 3}

\subsection{Initial Reductions}

     Let $\Phi(t), W(t)$ be two real-valued and non-negative smooth functions compactly supported in $(0,1)$, satisfying $\Phi(t)=W(t)=1$ for $t \in (1/U, 1-1/U)$ and such that
    $\Phi^{(j)}(t), W^{(j)}(t)\ll_j U^j$ for all integers $j \geq 0$. We consider the following smoothed sum
\begin{align*}
   S_i(X,Y;U) =\sum_{n}\sum_{\chi \in S_{i,n}} \sum_{m} \chi(m) \Phi\left( \frac {n}{Y} \right) W \left( \frac {m}{X} \right), \quad i=3,4.
\end{align*}
   where $U$ is a parameter to be chosen later. Applying the P\'olya-Vinogradov inequality in a way similar to the argument described in the Introduction, we see that when $Y \leq X$, for any $\varepsilon>0$,
\begin{align}
\label{1stredn}
   \Big |S_i(X,Y)- S_i(X,Y;U) \Big | \ll \frac {XY^{1/i+\varepsilon}+Y^{3/2+\varepsilon}}{U}.
\end{align}

   As the treatments for $S_3(X,Y)$ and $S_4(X,Y)$ are similar, we shall concentrate on the proof of the case $S_3(X,Y)$ in what follows and discuss briefly the proof of $S_4(X,Y)$ at the end of the paper.  Let $\chi$ be a Dirichlet character of modulus $q$.  Then we have the following well-known Poisson summation formula:
\begin{equation}
\label{eq:Poisson1dim}
 \sum_{m \in \mz}\chi(m) w\leg{m}{M}  = \frac{M}{q} \sum_{k \in \mz}  \tau(k, \chi) \widetilde{w}\leg{kM}{q},
\end{equation}
where $\tau(k, \chi)$ is the Gauss sum defined in \eqref{taur} and $\widetilde{w}$ is the Fourier transform of $w$, that is,
\[ \widetilde{w}(x) = \int\limits_{-\infty}^{\infty}  w(y) e(-xy) \dif y. \]

    Applying \eqref{eq:Poisson1dim} and \eqref{taug3}, we see that
\begin{align*}
  S_3(X,Y;U) &=X \sum_{k \in
   \mz } \sum_{n}\sum_{\chi \in S_{3,n}} \frac {  \tau(k, \chi)}{n} \Phi \left( \frac {n}{Y} \right)\widetilde{W} \leg{kX}{n} \\
   & =X \sum_{k \in
   \mz } \ \sumprime_{n \equiv 1 \bmod 3} \frac { \overline{\leg {\sqrt{-3}}{n}_3}g_3(k, n)}{N(n)} \Phi \left( \frac {N(n)}{Y} \right)\widetilde{W} \left(\frac {kX}{N(n)}\right)  =: M_0+R,
\end{align*}
   where $\Sigma'$ indicates that the sum is over $n \in \mz[\omega]$ with no $\mq$-rational prime divisor,
   \[
  M_0 =X\widetilde{W} \left(0\right) \sumprime_{n \equiv 1 \pmod 3} \frac { \overline{\leg {\sqrt{-3}}{n}}_3 g_3(0, n)}{N(n)} \Phi \left( \frac {N(n)}{Y} \right), \]
  and
\[  R =X \sum_{\substack{k \in
   \mz \\ k \neq 0} } \sumprime_{n \equiv 1 \pmod 3} \frac { \overline{\leg {\sqrt{-3}}{n}}_3 g_3(k, n)}{N(n)} \Phi \left( \frac {N(n)}{Y} \right)\widetilde{W} \left(\frac {kX}{N(n)}\right). \]

\subsection{The Term $M_{0}$}

    We estimate $M_{0}$ first. It follows directly from the definition that $g_3(0,n)=\varphi_{\omega}(n)$ if $n$ is a cubic power and $g_3(0,n)=0$ otherwise.
    Here $\varphi_{\omega}(n)$ denotes the number of reduced residue classes in $\mz[\omega]/(n)$.  Thus
\begin{align*}
  M_{0}= X\widetilde{W}(0)\sumprime_{\substack {n \equiv 1 \bmod {3} \\  \text{$n$ a cubic}}}\frac {\varphi_{\omega}(n)}{N(n)}\Phi \left( \frac {N(n)}{Y} \right).
\end{align*}
  As it is easy to see that $n^3$ has no $\mq$-rational prime divisor if and only if $n$ has no $\mq$-rational prime divisor, we can recast $M_0$ upon replacing $n$ by $n^3$ as
\begin{align*}
  M_{0}= X\widetilde{W}(0)\sumprime_{\substack {n \equiv 1 \bmod {3}}}\frac {\varphi_{\omega}(n^3)}{N(n^3)}\Phi \left( \frac {N(n^3)}{Y} \right)
  =X\widetilde{W}(0)\sumprime_{\substack {n \equiv 1 \bmod {3}}}\frac {\varphi_{\omega}(n)}{N(n)}\Phi \left( \frac {N^3(n)}{Y} \right).
\end{align*}

    By Mellin inversion, we have
\begin{align*}
    \Phi \left( \frac {N^3(n)}{Y} \right) = \frac 1{2\pi i}\int\limits_{(2)} \left( \frac{Y}{N^3(n)} \right)^s\widehat{\Phi}(s) \dif s \quad \mbox{where} \quad   \widehat{\Phi}(s)=\int\limits^{\infty}_{0}\Phi(t)t^{s-1} \dif t.
\end{align*}

     Integration by parts shows $\widehat{\Phi}(s)$ is a function satisfying the bound for all $\Re(s) > 0$, and integers $A>0$,
\begin{align}
\label{boundsforphi}
  \widehat{\Phi}(s) \ll (1+|s|)^{-A} U^{A-1}.
\end{align}
We then deduce that
\begin{align}
\label{M0}
  M_{0}=  \frac{X\widetilde{W}(0)}{2\pi i}\int\limits_{(2)}Y^s\widehat{\Phi}(s)\left( \sumprime_{n \equiv 1 \bmod 3}\frac {\varphi_{\omega}(n)}{N(n)^{1+3s}}
   \right) \dif s .
\end{align}

   Note that (see \cite[Prop. 9.1.2]{I&R}) a prime $\varpi \in \mz[\omega]$ is not $\mq$-rational if and only if $N(\varpi)=p \equiv 1 \pmod 3$ . We then deduce that when $\Re(s)>1$,
\begin{align*}
 & \sumprime_{n \equiv 1 \bmod 3}\frac {\varphi_{\omega}(n)}{N(n)^{1+s}} =  \sum_{\substack{ n \\ \varpi | n \Rightarrow N(\varpi)=p \equiv 1 \bmod 3}}
 \frac {\varphi_{\omega}(n)}{N(n)^{1+s}} \\
 = & \prod_{\substack{p \\ p \equiv 1 \bmod 3}} \left( 1+ \left( 1-\frac 1p \right)\frac {p^{-s}}{1-p^{-s}} \right)
 =\prod_{\substack{p}}\left(1+\left ( \frac {\chi_{0,3}(p)+\chi_{1,3}(p)}{2} \right) \left( 1-\frac 1p \right) \frac {p^{-s}}{1-p^{-s}} \right ),
\end{align*}
   where $\chi_{0,3}$ is the principal Dirichlet character modulo $3$ and $\chi_{1,3}$ the non-principal Dirichlet character modulo $3$.  Let $L(s, \chi_{i,3})$ stand for the corresponding Dirichlet $L$-functions for $i=0$, $1$. Now we define for $\Re(s)>1$,
\begin{align*}
   f(s)= L^{-1}(s, \chi_{0,3})L^{-1}(s, \chi_{1,3}) \prod_{\substack{p}}\left(1+\left ( \frac {\chi_{0,3}(p)+\chi_{1,3}(p)}{2} \right )
  \left( 1-\frac 1p \right)\frac {p^{-s}}{1-p^{-s}} \right )^2 .
\end{align*}
Observe that $f(x)=\prod_pf_p(s)$ with
\begin{align*}
  f_p(s)= &  \left(1+\frac {\chi_{0,3}(p)+\chi_{1,3}(p)}{p^{2s}}\cdot \frac {1-p^{s-1}}{1-p^{-s}}-
 \frac {(\chi_{0,3}(p)+\chi_{1,3}(p))^2}{4} \left(1-\frac 1p \right)
 \frac {1}{p^{2s}(1-p^{-s})}  \right. \\
 & \hspace*{1cm} +\frac {\chi_{0,3}(p)\chi_{1,3}(p)}{p^{2s}}+\frac {(\chi_{0,3}(p)+\chi_{1,3}(p))\chi_{0,3}(p)\chi_{1,3}(p)}{4} \left( 1-\frac 1p \right)
 \frac {1}{p^{3s}(1-p^{-s})} \\
 & \hspace*{1cm} \left. +\frac {(\chi_{0,3}(p)+\chi_{1,3}(p))^2}{4}\left( 1-\frac 1p \right)
 \frac {1}{p^{2s} (1-p^{-s} )^2}\left(1-\frac {\chi_{0,3}(p)+\chi_{1,3}(p)}{p^s}+\frac {\chi_{0,3}(p)\chi_{1,3}(p)}{p^{2s}} \right ) \right).
\end{align*}
  It follows from the expression of $f_p(s)$ that $f(s)$ is analytic for $\Re(s)>1/2$. \newline

   We then derive from \eqref{M0} that
\begin{align}
\label{M0g}
  M_{0}
  = X\widetilde{W}(0)\frac 1{2\pi i}\int\limits_{(2)}\frac {Y^sg(s)}{\sqrt{s-1/3}} \dif s,
\end{align}
   where
\begin{align*}
  g(s)=\widehat{\Phi}(s)\sqrt{ \left( s-\frac 13 \right)L(3s, \chi_{0,3})L(3s, \chi_{1,3})f(3s)}.
\end{align*}

   It is easy to see that $g(s)$ is analytic in a neighbourhood of $1/3$ and that
   \begin{equation} \label{gform}
   g(s)=g(1/3)+O(|s-1/3|)
   \end{equation}
    when $s$ is near $1/3$. \newline

    We now follow a method of E. Landau \cite{Landau1} (see also \cite[p. 187, exercise 21]{MVa1})
    by shifting the contour of integration in \eqref{M0g} to $\mathcal{C'} \cup \mathcal{C}$,
where $\mathcal{C'}$ is the contour running from $1/3-\varepsilon_0-i\infty$ to $1/3 -\varepsilon_0 - i\delta$ vertically and from $1/3-\varepsilon_0 + i\delta$ to
$1/3-\varepsilon_0+i\infty$ vertically. The contour $\mathcal{C}$ is from $1/3-\epsilon_0-i\delta$ to $1/3 - i\delta$ horizontally,
then along the semicircle $1/3 + \delta e(i\theta)$ , $-1/4 \leq \theta \leq 1/4$, and finally along a horizontal line to $1/3 -\varepsilon_0+ i\delta$.
Here $\varepsilon_0$ is sufficiently small and $\delta = 1/ \log Y$. \newline

    To estimate the integral over $\mathcal{C'}$,  we note that $L(s, \chi_{1,3})$ is bounded when $\Re(s)>0$ and that (see \cite[p. 100, exercise 3, ]{iwakow})
    when $\sigma=\Re(s) \geq 0$,
\begin{align}
\label{Lchibound}
  L(s, \chi_{i,3}) \ll (1+|s|)^{(1-\sigma)/2+\epsilon}, \quad i=0,1.
\end{align}

    We now divide the integral over $\mathcal{C'}$ into two parts,
\[ \int\limits_{\mathcal{C}'} = \int\limits_{\mathcal{C}'_1} + \int\limits_{\mathcal{C}'_2} ,\]
where $\mathcal{C}'_1$ is the part of $\mathcal{C'}$ with $|\Im(s)| \leq T$ ($T$ is to be chosen later) and $\mathcal{C}'_2$ being the rest. Applying \eqref{Lchibound} with \eqref{boundsforphi} by taking $A=1$ for the integral over $\mathcal{C}'_1$ and
\eqref{Lchibound} with \eqref{boundsforphi} by taking $A=2$ for the integral over $\mathcal{C}'_2$, we deduce that the integral over $\mathcal{C'}$ is
\begin{align} \label{estC}
    \ll XY^{1/3-\epsilon_0} \left( T^{2\epsilon_0}+\frac {U}{T^{1-2\varepsilon_0}} \right) \ll XY^{1/3-\varepsilon_0}U^{2\varepsilon_0} ,
\end{align}
upon setting $T=U$. \newline

Next, the integral over $\mathcal{C}$ is, using \eqref{gform},
\[   =\frac 1{2\pi i}\int\limits_{\mathcal{C}}\frac {Y^sg(1/3)}{\sqrt{s-1/3}} \dif s+O \left( \int\limits_{\mathcal{C}}\frac {Y^{\sigma}|s-1/3|}{\sqrt{s-1/3}} \dif |s| \right). \]
After two changes of variables, we get that
\begin{equation} \label{hankint}
\frac 1{2\pi i}\int\limits_{\mathcal{C}}\frac {Y^sg(1/3)}{\sqrt{s-1/3}} \dif s = \frac {Y^{1/3} g(1/3)}{2\pi i \sqrt{\log Y} }\int\limits_{\mathcal{D}} e^w w^{-1/2} \dif w
\end{equation}
where $D$ is the contour that goes from $-\varepsilon_0\log Y-i$ to $-i$ horizontally, then to $i$ along the right half of the unit circle centered at the origin and then to $-\varepsilon_0\log Y+i$ horizontally.  Using \cite[Theorem C.3]{MVa1} with $s = 1/2$, the right-hand side of \eqref{hankint} is
\[ \frac {g(1/3)Y^{1/3}}{\sqrt{ \log Y}} \left(  \frac{1}{\Gamma(1/2) } - \frac {1}{2\pi i  }\int\limits_{\mathcal{D}_+ \cup \mathcal{D}_-} e^w w^{-1/2} \dif w \right), \]
where $\mathcal{D}_{\pm}$ are the horizontal rays from $-\varepsilon_0 \log Y \pm i$ to $- \infty \pm i$ with the appropriate directions.  The integrals over $\mathcal{D}_{\pm}$ can be easily estimated as $\ll Y^{-\varepsilon_0}$.  Thus
\begin{align}
\label{estC'1}
    \frac 1{2\pi i}\int\limits_{\mathcal{C}}\frac {Y^sg(1/3)}{\sqrt{s-1/3}} \dif s=\frac {g(1/3)Y^{1/3}}{\sqrt{\pi \log Y}}+O\left( Y^{1/3-\varepsilon_0} \right).
\end{align}

   On the other hand, it is easy to show that
\begin{align}
\label{estC'2}
   \int\limits_{\mathcal{C}}\frac {Y^{\sigma}|s-1/3|}{\sqrt{s-1/3}} \dif |s| \ll \frac {Y^{1/3}}{(\log Y)^{3/2}} .
\end{align}

    Combining \eqref{estC}, \eqref{estC'1} and \eqref{estC'2}, we see that
\begin{align*}
  M_{0}
  = g \left(\frac 13 \right)\widetilde{W}(0)\frac {XY^{1/3}}{\sqrt{\pi \log Y}}+O \left( \frac {XY^{1/3}}{(\log Y)^{3/2}}+XY^{1/3-\varepsilon_0}U^{2\varepsilon_0} \right).
\end{align*}

   Using the observations that
\begin{align*}
  \widehat{\Phi} \left( \frac 13 \right)=3+O \left(\frac{1}{U^{1/3}} \right) \quad \mbox{and} \quad \widetilde{W}(0)=1+O\left( \frac 1U \right),
\end{align*}
   we can further rewrite $M_0$ as
\begin{align}
\label{M0asymp}
  M_{0}
  = C_1\frac {XY^{1/3}}{\sqrt{\log Y}}+O \left( \frac {XY^{1/3}}{(\log Y)^{3/2}}+\frac {XY^{1/3}}{(\log Y)^{1/2}U^{1/3}}+XY^{1/3-\epsilon_0}U^{2\varepsilon_0} \right),
\end{align}
   where
\begin{align}
\label{C1}
  C_1=\frac {3}{\sqrt{\pi}}\sqrt{\lim_{s \rightarrow 1/3} \left( s-\frac 13 \right)L(3s, \chi_{0,3})L(3s, \chi_{1,3})f(3s)}.
\end{align}
   We note that as $\chi_{1,3}$ is quadratic, it is easy to see that $C_1>0$.

\subsection{The Term $R$}
    Now suppose $k \neq 0$. We first apply the M\"obius function to detect the condition that $n \equiv 1 \pmod 3$ has no rational prime divisor to see that
\begin{align*}
  R & =X \sum_{\substack{k \in
   \mz \\ k \neq 0} } \  \sum_{\substack{ d \in \mz \\ d \equiv 1 \bmod 3}} \mu_{\mz}(d) \sum_{\substack{ n \in \mz[\omega] \\ n \equiv 1 \bmod 3}} \frac { \overline{\leg {\sqrt{-3}}{nd}}_3 g_3(k, nd)}{N(nd)} \Phi \left( \frac {N(nd)}{Y} \right)\widetilde{W} \left(\frac {kX}{N(nd)}\right).
\end{align*}
   where we define $\mu_{\mz}(d)=\mu(|d|)$, the usual M\"obius function. \newline

  We now apply \eqref{gprod} to conclude that
\begin{align*}
R &= X \sum_{\substack {k \in
   \mz \\ k \neq 0}} \sum_{\substack{ d \in \mz \\ d \equiv 1 \bmod 3}} \mu_{\mz}(d)\frac { \overline{\leg {\sqrt{-3}}{d}}_3 g_3(k, d)}{N(d)} H(k,d; X,Y),
\end{align*}
  where
\begin{align*}
H(k,d; X, Y)= \sum_{\substack{ n \in \mz[\omega] \\ n \equiv 1 \bmod 3}} \frac { \overline{\leg {d\sqrt{-3}}{n}}_3 g_3(k, n)}{N(n)} \Phi \left( \frac {N(nd)}{Y} \right)\widetilde{W} \left(\frac {kX}{N(nd)}\right).
\end{align*}

   We now write
\begin{align*}
 R = R_1(Z)+R_2(Z),
\end{align*}
   with
\[ R_1(Z) = X \sum_{\substack {k \in \mz \\ k \neq 0}} \ \sum_{\substack{ d \in \mz, \ |d| \leq Z\\ d \equiv 1 \bmod 3}} \mu_{\mz}(d)\frac { \overline{\leg {\sqrt{-3}}{d}}_3 g_3(k, d)}{N(d)} H(k,d; X,Y) \]
   and
   \[
R_2(Z) = X \sum_{\substack {k \in
   \mz \\ k \neq 0}} \ \sum_{\substack{ d \in \mz, \ |d| > Z \\ d \equiv 1 \bmod 3}} \mu_{\mz}(d)\frac { \overline{\leg {\sqrt{-3}}{d}}_3 g_3(k, d)}{N(d)} H(k,d; X,Y) .
\]

    We first estimate $R_2(Z)$ by noting that it follows from the definition of $\Phi$ that $H=0$ unless $N(nd) \leq Y$.
    We also note that if $d$ is square-free as a rational integer, it is also square-free in $\mz[\omega]$.  Hence it follows from \eqref{grnbound} that
\begin{align}
\label{gdbound}
    g_3(k,d) \leq N(d)^{1/2}.
\end{align}

    We further note that it follows from the definition of $W$ and integration by parts that for any $l \geq 0$, $j \geq 1$,
\begin{align}
\label{Wbound}
   \widetilde{W}^{(l)}(t) \ll \frac {U^{j-1}}{|t|^j}.
\end{align}

    Applying the above bound with $l=0$, $j=1$ or $j=2$ together with the trivial bound $g_3(k,n) \leq N(n)$,  we deduce that
\begin{align*}
H(k,d; X, Y) \ll & \sum_{N(n) \leq  Y/d^2} \left| \widetilde{W} \left(\frac {kX}{N(nd)}\right) \right| \ll \min \left(\sum_{N(n) \leq  Y/d^2}\frac {N(nd)}{|k|X},
\sum_{N(n) \leq  Y/d^2}\left(\frac {N(nd)}{kX}\right)^2U \right )
\\ \ll & \min \left(\sum_{N(n) \leq  Y/d^2}\frac {Y}{|k|X}, \sum_{N(n) \leq  Y/d^2}\left(\frac {Y}{kX}\right)^2U \right )\ll
\min \left(\frac {Y^2}{|k|d^2X}, \frac {Y}{d^2}\left(\frac {Y}{kX}\right)^2U \right ).
\end{align*}

Using \eqref{gdbound} and the above bound for $H$ leads us to the bound
\begin{align}
\label{R2}
R_2(Z) \ll X \sum_{\substack {k \in
   \mz \\ k \neq 0 \\ |k| \leq U}} \sum_{\substack{ d \in \mz \\ d \equiv 1 \bmod 3 \\ |d| > Z}} \frac {Y^2}{|k|d^3X}+ X \sum_{\substack {k \in
   \mz \\ k \neq 0 \\ |k|>U}} \sum_{\substack{ d \in \mz \\ d \equiv 1 \bmod 3 \\ |d| > Z}}  \frac {Y}{d^3}\left(\frac {Y}{kX}\right)^2U
   \ll  \frac {XY}{Z^2}\left(\frac {Y}{X}\right)^2U^{\varepsilon}.
\end{align}

Now, to estimate $R_1(Z)$,  we apply the Mellin inversion to obtain
\[ \Phi \left( \frac {N(n)}{Y} \right) \widetilde{W}\left(\frac {kX}{N(n)}\right) = \frac 1{2\pi i}\int\limits_{(2)} \left( \frac{Y}{N(n)} \right)^s\tilde{f}(s,k) \dif s, \quad \mbox{where} \quad  \tilde{f}(s,k)=\int\limits^{\infty}_{0}\Phi(t)\widetilde{W}\left(\frac {kX}{Yt}\right) t^{s-1} \dif t. \]

     Integration by parts and using \eqref{Wbound} shows $\tilde{f}(s)$ is a function satisfying the bound
\begin{align}
\label{boundsforf}
  \tilde{f}(s,k) \ll (1+|s|)^{-D} \left( 1+\frac {|k|X}{Y} \right)^{-E+D} U^{E-1},
\end{align}
for all $\Re(s) > 0$, and integers $D \geq 0, E>0$.   We deduce from the above discussions and \eqref{eq:gmult} that
\begin{align*}
   H(k,d; X, Y) &=  \frac 1{2\pi i}\int\limits_{(2)}\tilde{f}(s,k)\left ( \frac {Y}{d^2} \right )^sG(1+s,dk) \dif s, \quad \mbox{where} \quad
  G(1+s,dk) = \sum_{\substack{n \in \mz[\omega], \ (n,d)=1 \\ n \equiv 1 \bmod 3}} \frac { \overline{\leg {\sqrt{-3}}{n}}_3 g_3(dk, n)}{N(n)^{1+s}}.
\end{align*}

     We now move the line of integration to the line $\Re(s) = \varepsilon$. By Lemma \ref{lem1}, the only possible poles are at $s = 1/3$. Thus we may write $R_1(Z) = R_{1,1}(Z)+R_{1,2}(Z)$, where
\[  R_{1,1}(Z) =  X Y^{1/3} \sum_{\substack {k \in
   \mz \\ k \neq 0}} \ \sum_{\substack{ d \in \mz, \ |d| \leq Z \\ d \equiv 1 \bmod 3}} \mu_{\mz}(d)\frac { \overline{\leg {\sqrt{-3}}{d}}_3 g_3(k, d)}{N(d)d^{2/3}} \tilde{f} \left( \frac1{3}, k \right)\text{Res}_{s=1/3} G(1+s, dk), \]
   and
\[   R_{1,2}(Z) = \frac{X}{2\pi i} \sum_{\substack {k \in
   \mz \\ k \neq 0}} \ \sum_{\substack{ d \in \mz, \ |d| \leq Z \\ d \equiv 1 \bmod 3}} \mu_{\mz}(d)\frac { \overline{\leg {\sqrt{-3}}{d}}_3 g_3(k, d)}{N(d)}  \int\limits_{(\varepsilon)}\tilde{f}(s,k)\left ( \frac {Y}{d^2} \right )^sG(1+s,dk) \dif s. \]

     To estimate $R_{1,1}(Z)$, we note first that by the remark above \eqref{grnbound}, we may restrict the sum over $d$ to those satisfying $(d,k)=1$.
We then apply Lemma \ref{lem1} and the bound \eqref{boundsforf} with $D=0, E=1$ to obtain
\begin{align}
\label{R11}
   R_{1,1}(Z) \ll &  X Y^{1/3} \sum_{\substack {k=k_1k_2^2k_3^3 \in
   \mz \\ k \neq 0}} \left ( \frac {Y}{|k|X} \right ) \sum_{\substack{ d \in \mz, \ (d,k)=1 \\ 0 \neq d \leq Z}} \frac {N(dk_1)^{-1/6+\varepsilon}}{d^{5/3}}
   \ll  X Y^{1/3}\left ( \frac {Y}{X} \right ).
\end{align}

    To estimate $R_{1,2}(Z)$, we apply Lemma \ref{lem1} and bound \eqref{boundsforf} with $D=2$, $E=4$ to get that
\begin{equation} \label{R12}
   R_{1,2}(Z) \ll X Y^{\varepsilon} \sum_{\substack {k \in
   \mz \\ k \neq 0 }}\left ( \frac {Y}{kX} \right )^2U^3 \sum_{\substack{ d \in \mz \\ 0 \neq |d| \leq Z}} \frac {N(d)^{1/4}N(k)^{1/4}}{|d|}
   \int\limits_{\mr}(1+|t|)^{-1-\varepsilon} \dif t  \ll X Y^{\varepsilon} \left ( \frac {Y}{X} \right )^2 U^{3}Z^{1/2}.
\end{equation}

   Combining \eqref{R2}, \eqref{R11} and \eqref{R12}, we conclude that
\begin{align}
\label{R}
   R \ll &  \frac {XY}{Z^2}\left(\frac {Y}{X}\right)^2U^{\varepsilon}+X Y^{1/3}\left ( \frac {Y}{X} \right )+
   X Y^{\varepsilon} \left ( \frac {Y}{X} \right )^2 U^{3}Z^{1/2}.
\end{align}

\subsection{Conclusion}
   We now combine \eqref{1stredn}, \eqref{M0asymp} and \eqref{R} and adjust the value of $\varepsilon_0$ to see that
\begin{equation} \label{S3final}
\begin{split}
  S_3(X,Y)= C_1\frac {XY^{1/3}}{\sqrt{\log Y}}+& O\left( \frac {XY^{1/3}}{(\log Y)^{3/2}}+\frac {XY^{1/3}}{(\log Y)^{1/2}U^{1/3}}+ X Y^{1/3} \left( \frac {Y}{X}  \right) + XY^{1/3-\varepsilon}U^{2\varepsilon}  \right.\\
  & \hspace*{1cm} \left. +  \frac {XY^{1/3+\varepsilon}+Y^{3/2+\varepsilon}}{U}+\frac {XY}{Z^2}\left(\frac {Y}{X}\right)^2U^{\varepsilon}+
   X Y^{\varepsilon} \left ( \frac {Y}{X} \right )^2 U^{3}Z^{1/2} \right) ,
\end{split}
\end{equation}
where $C_1$ is given in \eqref{C1}.  As $Y \geq X^{6/7}$ implies that $Y^{3/2} \geq XY^{1/3}$, we now choose the values of $U, Z$ so that
\begin{align}
\label{UZ1}
  \frac {Y^{3/2}}{U}= \frac {XY}{Z^2}\left(\frac {Y}{X}\right)^2=
   X\left ( \frac {Y}{X} \right )^2 U^{3}Z^{1/2} .
\end{align}

   We then deduce that
\begin{align*}
   U = Y^{3/34}\left(\frac {X}{Y}\right)^{5/17}.
\end{align*}

  Substituting the above value of $U$ into \eqref{S3final},  making use of \eqref{UZ1} and our assumption that $Y \geq X^{6/7}$, we arrive at the
  expression given in \eqref{S3asymp}. \newline

We end this paper by giving a sketch on the proof of \eqref{S4asymp}. We proceed in a way similar to the proof \eqref{S3asymp} and the main term corresponding to $M_0$ given in \eqref{M0asymp} is
\begin{align*}
  C_2\frac {XY^{1/4}}{\sqrt{\log Y}}+O \left( \frac {XY^{1/4}}{(\log Y)^{3/2}}+\frac {XY^{1/4}}{(\log Y)^{1/2}U^{1/4}}+XY^{1/4-\epsilon_0}U^{3\varepsilon_0} \right),
\end{align*}
   where
\begin{equation}
\label{C2}
  C_2= \frac {4}{\sqrt{\pi}}\sqrt{\lim_{s \rightarrow 1/4} \left( s-\frac 14 \right)L(4s, \chi_{0,4})L(4s, \chi_{1,4})h(4s)},
  \end{equation}
  with
  \[   h(s)= L^{-1}(s, \chi_{0,4})L^{-1}(s, \chi_{1,4})\left ( \prod_{\substack{p}}\left(1+\left ( \frac {\chi_{0,4}(p)+\chi_{1,4}(p)}{2} \right )
  \left(1-\frac 1p \right)\frac {p^{-s}}{1-p^{-s}} \right ) \right )^2.  \]
   Here $\chi_{0,4}$ is the principal Dirichlet character modulo $4$ and $\chi_{1,4}$ the non-principal Dirichlet character modulo $4$.
$L(s, \chi_{i,4})$ is the corresponding Dirichlet $L$-functions for $i=0,1$.  It is also easy to see that $C_2>0$. \newline

   The estimation corresponding to $R_2$ for $S_3(X,Y)$ remains the same for $S_4(X,Y)$. The estimation corresponding to $R_{1,1}$ for $S_3(X,Y)$ in this case
becomes (note that by \eqref{g4prod} and \eqref{eq:gmult}, we have $g_4(d^2k, n)$ in place of $g_3(dk,n)$)
\begin{align*}
  X Y^{1/4} \sum_{\substack {k \in
   \mz, \ k \neq 0 \\ |k| \leq YU^3/X}} & \left ( \frac {Y}{|k|X} \right ) \sum_{\substack{ d \in \mz \\ 0 \neq |d| \leq Z}} \frac {N(d^2)^{1/8+\epsilon}N(k)^{1/8+\epsilon}}{|d|}+
   X Y^{1/4} \sum_{\substack {k \in
   \mz, \ k \neq 0 \\ |k| >YU^3/X}}\left ( \frac {Y}{|k|X} \right )^2U^3
   \sum_{\substack{ d \in \mz \\ 0 \neq |d| \leq Z}} \frac {N(d^2)^{1/8+\epsilon}N(k)^{1/8+\epsilon}}{|d|} \\
    \ll &
   X Y^{1/4} \left ( \frac {Y}{X} \right )^{5/4} U^{3/4+\epsilon}Z^{1/2+\epsilon},
\end{align*}
   provided that
\begin{align}
\label{Ucondition}
   YU^3/X \geq 1.
\end{align}

 The terms in $S_4(X,Y)$ analogous to $R_{1,2}$ in $S_3(X,Y)$ become (by applying Lemma \ref{lem1} and \eqref{boundsforf} with $D=3$, $E=5$)
\begin{align*}
   \ll  X Y^{\epsilon} \sum_{\substack {k \in
   \mz \\ k \neq 0}}\left ( \frac {Y}{kX} \right )^2U^4 \sum_{\substack{ d \in \mz \\ 0 \neq |d| \leq Z}} \frac {N(d^2)^{1/4}N(k)^{1/4}}{|d|}
   \int\limits_{\mr}(1+|t|)^{-3/2} \dif t
    \ll
   X Y^{\epsilon} \left ( \frac {Y}{X} \right )^2 U^4 Z.  \nonumber
\end{align*}

     We then conclude that when $Y \geq X^{4/5}$,
\begin{align*}
  S_4(X,Y)=C_2 \frac {XY^{1/4}}{\sqrt{\log Y}} +& O \left( \frac {XY^{1/4}}{(\log Y)^{3/2}} + \frac {XY^{1/4}}{(\log Y)^{1/2}U^{1/4}}+ XY^{1/4-\varepsilon}U^{3\varepsilon} \right. \\
   & \hspace*{1cm} \left. +X Y^{1/4}\left ( \frac {Y}{X} \right )^{5/4}U^{3/4+\varepsilon}Z^{1/2+\varepsilon}+ \frac {Y^{3/2+\varepsilon}}{U}+\frac {XY}{Z^2}\left(\frac {Y}{X}\right)^2U^{\varepsilon}+ X Y^{\varepsilon} \left ( \frac {Y}{X} \right )^2 U^4Z \right) ,
\end{align*}
   where $C_2$ is given in \eqref{C2}. \newline

   We now choose $U, Z$ such that
\begin{align*}
  X Y^{1/4}\left ( \frac {Y}{X} \right )^{5/4}U^{3/4}Z^{1/2}= \frac {XY}{Z^2}\left(\frac {Y}{X}\right)^2=X  \left ( \frac {Y}{X} \right )^2 U^4Z.
\end{align*}

  This implies that
\begin{align*}
  Z=Y^{9/31}\left (\frac {Y}{X}\right )^{12/31}, \quad U=Y^{1/31}\left (\frac {X}{Y} \right )^{9/31}.
\end{align*}

  One checks that \eqref{Ucondition} is satisfied when $Y \geq X^{4/5}$ and this leads to the expression for $S_4(X,Y)$ in \eqref{S4asymp}.  \newline

\noindent{\bf Acknowledgments.} P. G. is supported in part by NSFC grant 11871082 and L. Z. by the FRG grant PS43707 and the Faculty Silverstar Award PS49334 at the University of New South Wales (UNSW).  Parts of this work were done when P. G. visited UNSW in August 2018. He wishes to thank UNSW for the invitation, financial support and warm hospitality during his pleasant stay.  Both authors would like to thank the anonymous referee for his/her comments and suggestions.

\bibliography{biblio}
\bibliographystyle{amsxport}

\vspace*{.5cm}

\noindent\begin{tabular}{p{8cm}p{8cm}}
School of Mathematical Sciences & School of Mathematics and Statistics \\
Beihang University & University of New South Wales \\
Beijing 100191 China & Sydney NSW 2052 Australia \\
Email: {\tt penggao@buaa.edu.cn} & Email: {\tt l.zhao@unsw.edu.au} \\
\end{tabular}

\end{document}